\documentclass{IEEEtran4PSCC}
\ifCLASSINFOpdf
   \usepackage[pdftex]{graphicx}
\else
   \usepackage[dvips]{graphicx}
\fi
\usepackage[cmex10]{amsmath}

\usepackage{color}
\usepackage{amsmath,amsfonts}
\usepackage{algorithmic}
\usepackage{algorithm}
\usepackage{amssymb}
\usepackage{caption}
\usepackage{cite}
\allowdisplaybreaks

\newcommand{\calB}{\mathcal{B}}
\newcommand{\calE}{\mathcal{E}}
\newcommand{\calG}{\mathcal{G}}

\newcommand{\calI}{\mathcal{I}}
\newcommand{\calL}{\mathcal{L}}
\newcommand{\bd}{\mathbf{d}}

\makeatletter
\let\old@ps@headings\ps@headings
\let\old@ps@IEEEtitlepagestyle\ps@IEEEtitlepagestyle
\def\psccfooter#1{%
    \def\ps@headings{%
        \old@ps@headings%
        \def\@oddfoot{\strut\hfill#1\hfill\strut}%
        \def\@evenfoot{\strut\hfill#1\hfill\strut}%
    }%
    \def\ps@IEEEtitlepagestyle{%
        \old@ps@IEEEtitlepagestyle%
        \def\@oddfoot{\strut\hfill#1\hfill\strut}%
        \def\@evenfoot{\strut\hfill#1\hfill\strut}%
    }%
    \ps@headings%
}
\makeatother

\begin{document}
\title{GPU-Accelerated Sequential Quadratic Programming Algorithm for Solving ACOPF}

%% To specify the authors when (number of affiliations <= 2)
% \author{
% \IEEEauthorblockN{Author n.1 Name per Affiliation A\\ Author n.2 Name per Affiliation A}
% \IEEEauthorblockA{(Affiliation A) Department Name of Organization \\
% Name of the organization, acronyms acceptable\\
% City, Country\\
% \{email author n.1, email author n.2\}@domain (if desired)}
% \and
% \IEEEauthorblockN{Author n.1 Name per Affiliation B\\ Author n.2 Name per Affiliation B}
% \IEEEauthorblockA{(Affiliation B) Department Name of Organization \\
% Name of the organization, acronyms acceptable\\
% City, Country\\
% \{email author n.1, email author n.2\}@domain (if desired)}
% }
\author{
% \IEEEauthorblockN{Bowen Li}
% \IEEEauthorblockA{Imperial College Business School\\
% Imperial College London, London, UK\\
% libowen@umich.edu}
% \and
\IEEEauthorblockN{Bowen Li, \ Michel Schanen, \ Kibaek Kim}
\IEEEauthorblockA{Mathematics and Computer Science Division\\
Argonne National Laboratory, Lemont, IL, USA\\
libowen@umich.edu, \ \{mschanen, kimk\}@anl.gov}
}

%% To specify the authors when (number of affiliations > 2)
% \author{\IEEEauthorblockN{Author n.1\IEEEauthorrefmark{1},
% Author n.2\IEEEauthorrefmark{2},
% Author n.3\IEEEauthorrefmark{3}, 
% Author n.4\IEEEauthorrefmark{3} and
% Author n.5\IEEEauthorrefmark{4}}
% \IEEEauthorblockA{\IEEEauthorrefmark{1} Department Name of Organization A\\
% Name of the organization A,
% Address A\\ Emails if wanted}
% \IEEEauthorblockA{\IEEEauthorrefmark{2} Department Name of Organization B\\
% Name of the organization B,
% Address B\\ Emails if wanted}
% \IEEEauthorblockA{\IEEEauthorrefmark{3} Department Name of Organization C\\
% Name of the organization C,
% Address C\\ Emails if wanted}
% \IEEEauthorblockA{\IEEEauthorrefmark{4}Department Name of Organization D\\
% Name of the organization D,
% Address D\\ Emails if wanted}
% }

\maketitle

\begin{abstract}
Sequential quadratic programming (SQP)  is widely used in solving nonlinear optimization problem, with advantages of warm-starting solutions, as well as finding high-accurate solution and converging quadratically using second-order information, such as the Hessian matrix. In this study we develop a scalable SQP algorithm for solving the alternate current optimal power flow problem (ACOPF), leveraging the parallel computing capabilities of graphics processing units (GPUs). Our methodology incorporates the alternating direction method of multipliers (ADMM) to initialize and decompose the quadratic programming  subproblems within each SQP iteration into independent small subproblems for each electric grid component. We have implemented the proposed SQP algorithm using our portable Julia package ExaAdmm.jl, which solves the ADMM subproblems in parallel on all major GPU architectures. For numerical experiments, we compared three solution approaches: (i) the SQP algorithm with a GPU-based ADMM subproblem solver, (ii) a CPU-based ADMM solver, and (iii)   the QP solver Ipopt (the state-of-the-art interior point solver) and observed that for larger instances our GPU-based SQP solver efficiently leverages the GPU many-core architecture, dramatically reducing the solution time.
\end{abstract}


% Use this to place sponsorships
\thanksto{\noindent This research was supported by the Exascale Computing Project (17-SC-20-SC), a collaborative effort of the U.S. Department of Energy Office of Science and the National Nuclear Security Administration. This material is based upon work supported by the U.S. Department of Energy, Office of Science, under contract number DE-AC02-06CH11357.}

\section{Introduction}

Sequential quadratic programming (SQP) has been widely employed for solving general nonlinear optimization problems (NLPs) \cite{fletcher2002nonlinear, sqp_wright, SQP_boggs, nocedal1999numerical, conn2000trust}. By exploiting the second-order information (e.g., the Hessian matrix), this method iteratively solves quadratic programming (QP) subproblems in order to derive a good search direction to solution. SQP also offers distinct advantages including high solution accuracy and a quadratic convergence rate, as well as warm-starting solutions particularly for mixed-integer nonlinear programming (e.g., \cite{leyffer2001integrating,mahajan2021minotaur}). In practise, SQP can be used to solve NLP problems in different applications. One such application is the alternate optimal power flow problem (ACOPF) \cite{OPF_book}, which determines the optimal generation output to meet demand while satisfying the nonlinear physical laws of the grid.

Despite its distinct advantages, SQP can become computationally burdensome when dealing with large QPs, such as those encountered in ACOPF with large-scale networks. This can pose challenges in terms of time consumption and scalability. To resolve the issue, we propose to use decomposition and parallel computing to transform large QP subproblems into small and independent optimization problems that can be efficiently solved on graphics processing units (GPUs) in parallel. In \cite{mhanna2019Adaptive}, a component-decomposition method is proposed for ACOPF problems. In \cite{kim2021Leveraging,kim2022accelerated}, the authors implemented the component-decomposition method using a convergence-guaranteed computation procedure \cite{sun2021two} and parallel computing on GPUs. In this work we extend the notions of decomposition and parallel computation to the SQP algorithm structure, given its benefits on convergence and solution quality.

To summarize our contributions, we develop a GPU-accelerated SQP algorithm for ACOPF, whose QP subproblems are solved by a parallel and scalable algorithm. Specifically, instead of directly solving the QP, we employ the distributed algorithm known as the alternating direction method of multipliers (ADMM) \cite{ADMM_book} to decouple the constraints in the QP. By leveraging the idea of component decomposition \cite{mhanna2019Adaptive}, we break down the QP subproblem into a large number of small optimization subproblems for each electric grid component (e.g., generators, buses, and transmission lines). The subproblem solutions are independent, hence allowing for parallel computation. For the transmission line subproblems, we implement improvements to avoid ineffective solutions resulting from overconstraints. As the ADMM algorithm progresses, we eventually solve the QP subproblem in a distributed and parallel fashion. Furthermore, we implement the proposed SQP algorithm for both CPU and GPU architecture using our portable software package \texttt{ExaAdmm.jl} \cite{kim2022accelerated}. 

To demonstrate the effectiveness of the proposed SQP algorithm with ADMM QP solver (SQP-ADMM), we compare the performance of a centralized SQP algorithm with \texttt{Ipopt} \cite{ipopt} as the QP solver with that of the SQP-ADMM methods running on CPUs and GPUs, respectively. We test using ACOPF problems of varying grid sizes up to 6,468 buses and observe that as the network size increases, the GPU-based SQP-ADMM method exhibits greater computational advantages over centralized SQP methods and CPU-based SQP-ADMM methods. 

The rest of the paper is organized as follows. A preliminary introduction to the trust-region SQP algorithm is first given in Section~\ref{sec: sqp}.
In Section~\ref{sec: sqp_acopf} we introduce the ACOPF formulation and formulate the corresponding QP subproblem for ACOPF with modifications. 
In Section~\ref{sec: ADMM} we discuss the ADMM algorithm by component decomposition to solve the QP subproblem as well as the GPU implementation to achieve parallel computing. Section~\ref{sec: result} demonstrates our numerical results by comparing with SQP-ADMM methods on CPUs and GPUs with centralized SQP methods with Ipopt. We also discuss the limitations of the current work. Section~\ref{sec: con} briefly summarizes our conclusions and possible future work.

% The preliminary test results are summarized in Table~\ref{tab}.

% Paragraph of ACOPF \cite{OPF_book}:

% Paragraph of SQP: \cite{fletcher2002nonlinear,sqp_wright,SQP_boggs,nocedal1999numerical, conn2000trust}

% Paragraph of ADMM: \cite{mhanna2019Adaptive, ADMM_book}

% Paragraph of GPU: \cite{kim2021Leveraging,kim2022accelerated}

% To summarize our contribution, we developed 

\section{Preliminary: Trust-Region SQP}
\label{sec: sqp}

This section serves as a preliminary introduction to the trust-region SQP algorithm considered in this paper.
For the algorithm description, we consider the following NLP formulation:
\begin{subequations}
\label{eq: nlp}
\begin{align}
    \min_x \quad & f(x) \\
    \text{s.t.} \quad 
    & Ax \geq b, \label{eq: lin_con} \\
    & c_i(x) = 0,\quad \forall i\in\calE, \label{eq: nonlin_e}\\
    & c_i(x) \geq 0,\quad \forall i\in\calI, \label{eq: nonlin_i}
\end{align}
\end{subequations}
where $x\in\mathbb{R}^n$ is the decision variable with objective function $f:\mathbb{R}^n\rightarrow\mathbb{R}$. Equation \eqref{eq: lin_con} represents the linear constraints on $x$. Equations \eqref{eq: nonlin_e} and \eqref{eq: nonlin_i} respectively represent the nonlinear equality and inequality constraints on $x$ using the nonlinear function $c_i:\mathbb{R}^n\rightarrow\mathbb{R}$. Further, we assume all nonlinear functions $f$ and $c_i$ are twice continuously differentiable. 

The Lagrangian function of problem~\eqref{eq: nlp} is given by
\begin{align}
    \calL(x,y,\lambda) = f(x) + y^\top(Ax-b) + \sum_{i\in\calE\cup\calI} \lambda_i c_i(x),
\end{align}
where $y$ and $\lambda_i$ are the Lagrangian multipliers corresponding to linear constraints and nonlinear constraint $i$, respectively. Next, we give the key components of an SQP algorithm. 

\subsection{Linear Feasibility}
The linear feasibility problem projects the initial solution $x^{(0)}$ to the feasible set with respect to the linear constraints \eqref{eq: lin_con} by solving  
\begin{align} \label{eq: linfeas}
    \min \quad & \|x - x^{(0)}\|^2 \quad \text{s.t.} \quad Ax \geq b.
\end{align}
Once the solution of \eqref{eq: linfeas} is obtained, the subsequent steps of the SQP algorithm maintains the feasibility with respect to the linear constraints. Hence, the first QP subproblem is derived with this projection solution.
Note that the positive objective value of \eqref{eq: linfeas} gives the certification of the infeasibility with respect to the original problem \eqref{eq: nlp}.

\subsection{QP Subproblem}
Given solution $x_k$, the SQP algorithm approximates \eqref{eq: nlp} by the following QP subproblem: 
\begin{subequations}
\label{eq: qp}
\begin{align}
    \min_d \quad 
    & \nabla f_{k}^\top d + \frac{1}{2} d^\top \nabla_{xx}^2 \calL_{k} d \\
    \text{s.t.} \quad 
    & A (x^{(k)}+d) \geq b, \label{eq: lin_qp}\\
    & \nabla c_i(x^{(k)})^\top d + c_i(x^{(k)}) = 0, && \forall i\in\calE,\\
    & \nabla c_i(x^{(k)})^\top d + c_i(x^{(k)}) \geq 0, && \forall i\in\calI, \\
    & \|d\|_\infty \leq \Delta^{(k)}, \label{eq: TR}
\end{align}
\end{subequations}
where $d$ is a solution step from $x^{(k)}$, $\nabla f_{k} := \nabla f(x^{(k)})$ is the gradient of $f$ evaluated at $x^{(k)}$, $\nabla_{xx}^2\calL_{k}$ is the Hessian of the Lagrangian function $\calL$ evaluated at $(x^{(k)},y^{(k)},\lambda^{(k)})$, $\nabla c_i$ is the gradient of constraint $c_i$, and \eqref{eq: TR} is a $\ell_\infty$-norm trust-region constraint with radius $\Delta^{(k)}$. Note that the QP subproblem~\eqref{eq: qp} can be nonconvex when $\nabla_{xx}^2 \calL_{k}$ is not positive definite.

\subsection{Merit Functions}
The QP subproblem finds a trial step $d^{(k)}$ at the current iterate $x^{(k)}$, if feasible, such that $x^{(k+1)} = x^{(k)}+d^{(k)}$ may be the new trial iterate.
We use the $\ell_1$ merit function \cite{SQP_boggs, nocedal1999numerical}
\begin{align*}
    \phi_1(x;\mu) = f(x) + \mu h(x),
\end{align*}
where $h(x) := \sum_{i\in\calE} |c_i(x)| + \sum_{i\in\calI} [c_i(x)]^-$ and $[z]^-:= \max\{0,-z\}$. Note that $h(x)$  considers only nonlinear constraints in \eqref{eq: nlp} because linear constraints are guaranteed to be satisfied by the new trial iteration from \eqref{eq: lin_qp}.

The new trial iterate $x^{(k+1)}$ is accepted if $\phi_1(x^{(k)};\mu) - \phi_1(x^{(k)}+d^{(k)};\mu) > 0$ and
\begin{align*}
    \frac{\phi_1(x^{(k)};\mu) - \phi_1(x^{(k)}+d^{(k)};\mu)}{q(0; \mu) - q(d^{(k)}; \mu)} > 0,
\end{align*}
where $q(d;\mu)$ is a quadratic model of $\phi_1$ as given by
\begin{align*}
    q(d;\mu) &:= f_k + \nabla f_k^\top d + \frac{1}{2} d^\top \nabla_{xx}^2\calL_k d + \mu m(d) \\
    m(d) := &\sum_{i\in\calE} |c_i(x^{(k)})+\nabla c_i(x^{(k)})^\top d|
    \notag\\
    &+ \sum_{i\in\calI} \left[ c_i(x^{(k)})+\nabla c_i(x^{(k)})^\top d \right]^-.
\end{align*}
Otherwise, the new trial point is rejected.

\subsection{Algorithm Description}

\begin{algorithm}[ht!]
\caption{Basic trust-region SQP}
\begin{algorithmic}[1]\label{alg: sqp}
\STATE Given $x^{(0)}$ and $\Delta^{(0)}$, set $k=0$.
\IF{$x^{(0)}$ is infeasible to linear constraints}
    \STATE Solve linear feasibility problem~\eqref{eq: linfeas}.
\ENDIF
\REPEAT
    \STATE Evaluate $\nabla f_k$, $c_i(x^{(k)})$, $\nabla c_i(x^{(k)})$, and $\nabla_{xx}^2\calL_k$.
    \STATE Solve QP subproblem~\eqref{eq: qp} to find $d^{(k)}$.
    \IF{$d^{(k)}$ is accepted}
            \STATE Set $x^{(k+1)} \gets x^{(k)} + d^{(k)}$.
            \STATE Set $\Delta^{(k+1)} \geq \Delta^{(k)}$.
            % \IF{$\Delta^{(k)} = \|p^{(k)}\|_\infty$}
            %     \STATE Set $\Delta^{(k+1)} \gets \max\{2\Delta^{(k)},\bar\Delta\}$
            % \ENDIF
    \ELSE
            \STATE Set $\Delta^{(k+1)} \leq \Delta^{(k)}$.
                % \STATE Set $\Delta^{(k+1)} \gets 0.5 \min\{\Delta^{(k)},\|p^{(k)}\|_\infty\}$
    \ENDIF
    \STATE Set $k \gets k+1$.
\UNTIL{termination criteria satisfied}
\end{algorithmic}
\end{algorithm}

In Algorithm \ref{alg: sqp} we give the basic trust-region SQP based on the preceding sections.
Note that the algorithm  assumes that the linear feasibility problem \eqref{eq: linfeas} and the QP subproblem \eqref{eq: qp} are always solvable. As for termination criteria, we stop when the accepted step $d^{(k)}$ is smaller than tolerance or the optimality condition of \eqref{eq: acopf} is satisfied.
In this paper, the algorithm and our implementation focus on core functionality (i.e., solving the QP subproblem on GPU) for our numerical experiments, leaving as future work some advanced features such as second-order correction and feasibility restoration.

% In the case that problems \eqref{eq: linfeas} and the QP subproblem \eqref{eq: qp} are not solvable, Algorithm \ref{alg: sqp} will terminate early. 

\section{Implementation of SQP for ACOPF}\label{sec: sqp_acopf}

In this section we describe the ACOPF problem formulation in the context of the SQP algorithm described in Section \ref{sec: sqp}.
In particular, we focus on solving the QP subproblem \eqref{eq: qp} and present the QP subproblem formulation associated with the ACOPF problem. In addition, we show that, by eliminating some of the equality constraints, the solutions of the QP subproblem can be more effective in the progression of the SQP algorithm.

\subsection{ACOPF Formulation}
We present the ACOPF formulation used through the section.
ACOPF seeks the optimal planning of the generation units to minimize the cost while dispatching the power over the network with physics constraints. A rectangular formulation of ACOPF problem is given as
\begin{subequations}\label{eq: acopf}
\begin{align}
    \min \ 
    & \sum_{(g,i)\in\calG} f_{gi}(p_{gi}) \\
    \text{s.t.} \
    & \underline{p}_{gi} \leq p_{gi} \leq \bar{p}_{gi},\quad (g,i)\in\calG, \label{eq:opfc1}\\
    & \underline{q}_{gi} \leq q_{gi} \leq \bar{q}_{gi}, \quad (g,i)\in\calG, \label{eq:opfc2} \\
    & \underline{v}_i^2 \leq w_i \leq \bar{v}_i^2,\quad i\in\calB, \label{eq:opfc3}\\
    & -2\pi \leq \theta_i \leq 2\pi,\quad i\in\calB, \label{eq:opfc4}\\
    % & \tan(\underline{\theta}_{ij}^\Delta) w_{ij}^R \leq w_{ij}^I \leq \tan(\bar{\theta}_{ij}^\Delta) w_{ij}^R && (i,j)\in\calL \\
    % & \underline{\theta}_{ij}^\Delta \leq \theta_i - \theta_j \leq \bar{\theta}_{ij}^\Delta && (i,j)\in\calL \\
    & \sum_{(g,i)\in\calG} p_{gi} - p_i^d = \sum_{j\in\calB_i} p_{ij} + g_i^{sh} w_i, \quad i\in\calB, \label{eq:opfc5} \\
    & \sum_{(g,i)\in\calG} q_{gi} - q_i^d = \sum_{j\in\calB_i} q_{ij} - b_i^{sh} w_i, \quad i\in\calB, \label{eq:opfc6} \\
    % & p_{ij} = g_{ii} w_i + g_{ij} w_{ij}^R + b_{ij} w_{ij}^I && (i,j) \in \calL \\
    % & q_{ij} = - b_{ii} w_i - b_{ij} w_{ij}^R + g_{ij} w_{ij}^I && (i,j) \in \calL \\
    % & p_{ji} = g_{jj} w_j + g_{ji} w_{ij}^R - b_{ji} w_{ij}^I && (i,j) \in \calL \\
    % & q_{ji} = -b_{jj} w_j - b_{ji} w_{ij}^R - g_{ji} w_{ij}^I && (i,j) \in \calL \label{eq:opfc10}\\
    & (w_{ij}^R)^2 + (w_{ij}^I)^2 = w_i w_j, \quad  (i,j)\in\calL, \label{eq:opfc11}\\
    % & w_{ij}^R \sin_{ij} = w_{ij}^I \cos(\theta_i - \theta_j),\ (i,j)\in\calL \label{eq:opfc12} \\
    & w_{ij}^R \sin_{ij} = w_{ij}^I \cos_{ij},\ (i,j)\in\calL, \label{eq:opfc12} \\
    % & \theta_i - \theta_j = \atan2(w_{ij}^I,w_{ij}^R) && (i,j)\in\calL \label{eq:opfc12} \\
    & p_{ij}^2 + q_{ij}^2 \leq \bar{s}_{ij}^2,\quad (i,j)\in\calL, \label{eq:opfc13} \\
    & p_{ji}^2 + q_{ji}^2 \leq \bar{s}_{ji}^2,\quad (i,j)\in\calL, \label{eq:opfc14}
\end{align}
\end{subequations}
where the power flow variables are defined for the rectangular formulation as follows:
\begin{subequations}\label{eq: acopf_pf}
 \begin{align}
    & p_{ij} := g_{ii} w_i + g_{ij} w_{ij}^R + b_{ij} w_{ij}^I,\quad (i,j) \in \calL, \label{eq: pij_gb} \\
    & q_{ij} := - b_{ii} w_i - b_{ij} w_{ij}^R + g_{ij} w_{ij}^I,\quad (i,j) \in \calL, \label{eq: qij_gb}\\
    & p_{ji} := g_{jj} w_j + g_{ji} w_{ij}^R - b_{ji} w_{ij}^I,\quad (i,j) \in \calL, \label{eq: pji_gb}\\
    & q_{ji} := -b_{jj} w_j - b_{ji} w_{ij}^R - g_{ji} w_{ij}^I,\quad (i,j) \in \calL. \label{eq: qji_gb}
\end{align}   
\end{subequations}

We define the index set of transmission lines, buses, and generators as $\calL$, $\calB$, and $\calG$, respectively. 
In terms of decision variables, for the $g$th generator on bus $i$ denoted as $(g,i)\in\calG$, variables $p_{gi}$ and $q_{gi}$ represent its active and reactive power generation, which are bounded in \eqref{eq:opfc1} and \eqref{eq:opfc2}. We assume that the objective functiond $f_{gi}:\mathbb{R} \rightarrow \mathbb{R}$ are twice continuously differentiable. For each bus $i\in\calB$, variable $\theta_i$ represents the phase angle. Because of rectangular formulation~\cite{OPF_book}, we introduce variable $w_i$ to represent $v_i^2$, where $v_i$ is the voltage magnitude. Both magnitude and phase angle are bounded in \eqref{eq:opfc3} and \eqref{eq:opfc4}. Constraints \eqref{eq:opfc5} and \eqref{eq:opfc6} respectively represent the active and reactive power balance, where $p_{ij}$ and $p_{ji}$ are respectively the active and reactive power flow from bus $i$ to $j$ and set $B_i$ includes all nodes connect to bus $i$. For each transmission line $(i,j)\in\calL$, constraints \eqref{eq:opfc11} and \eqref{eq:opfc12} represent the relationship between phase angle $\theta$ and the rectangular-form variables $w$, where $\sin_{ij}:=\sin(\theta_i - \theta_j)$ and $\cos_{ij} := \cos(\theta_i - \theta_j)$. Variables $w_{ij}^R$ and $w_{ij}^I$ support the rectangular formulation and physically represent $v_iv_j\cos_{ij}$ and $v_iv_j\sin_{ij}$. Constraints \eqref{eq:opfc13} and \eqref{eq:opfc14} represent the line limit. 
Note that the only nonlinear constraints in \eqref{eq: acopf} and \eqref{eq: acopf_pf} are \eqref{eq:opfc11}--\eqref{eq:opfc14}. 

In terms of system parameters, we let phase angle at the reference bus $\theta_{ref} = 0$ as convention. For each generator $(g,i)\in\calG$, $\underline{p}_{gi}/\overline{p}_{gi}$ and $\underline{q}_{gi}/\overline{q}_{gi}$ are the lower and upper bounds of active and reactive power generation, respectively. For each bus $i\in\calB$, $\underline{v}_i$ and $\overline{v}_i$ respectively are the lower and upper bounds of the voltage magnitude; $p_i^d$ and $q_i^d$ respectively are the active and reactive power demand; and $g^{sh}_i + \mathbf{j}b^{sh}_i$ is the nodal shunt admittance. For each transmission line $(i,j)\in\calL$, $g_{ij}+\mathbf{j}b_{ij}$ represents the $i$th-row and $j$th-column entry in the admittance matrix $Y_{ij}$; $\bar{s}_{ij}^2$ and $\bar{s}_{ji}^2$ are the power flow limits.  

\subsection{ACOPF-Based QP Formulation}

We define the Lagrangian function  $\calL_{ij}$ evaluated at $(w_{ij}^{R}, w_{ij}^{I}, w_{i}, w_{j}, \theta_{i}, \theta_{j})$ with respect to constraints \eqref{eq:opfc11}--\eqref{eq:opfc14} and the corresponding multipliers $\pi$ for $(i,j)\in\calL$:
\begin{align}\label{eq: lag}
\calL_{ij} &:= \pi_{\ref{eq:opfc11}} \left[ (w_{ij}^R)^2 + (w_{ij}^I)^2 - w_i w_j \right] \nonumber\\
    & \quad + \pi_{\ref{eq:opfc12}} \left( w_{ij}^R \sin_{ij} - w_{ij}^I \cos_{ij} \right) \nonumber\\
    & \quad +  \pi_{\ref{eq:opfc13}} \left( \bar{s}_{ij}^2 - p_{ij}^2 - q_{ij}^2 \right) +  \pi_{\ref{eq:opfc14}} \left( \bar{s}_{ji}^2 - p_{ji}^2 - q_{ji}^2 \right),
\end{align}
where $p_{ij},q_{ij},p_{ji},q_{ji}$ are defined in \eqref{eq: acopf_pf}.
We denote by $\nabla^2\calL_{ij}^{(k)}$ the Hessian of the Lagrangian function at a given iterate $(p_{gi}^{(k)}, q_{gi}^{(k)}, w_{ij}^{R,(k)}, w_{ij}^{I,(k)}, w_{i}^{(k)}, \theta_{i}^{(k)})$ and the multiplier $\pi^{(k)}$.
Then, with the iterate, the QP subproblem of \eqref{eq: acopf} with trust-region $\Delta^{(k)}$ is given by 
\begin{subequations}\label{eq: QPsub}
\begin{align}
    & \min_{\bd,d_{w_{ij}^R},d_{w_{ij}^I}} \
    \sum_{(g,i)\in\calG} \left(\nabla f_{gi}(p_{gi}^{(k)})d_{p_{gi}} + \frac{1}{2}\nabla^2 f_{gi}(p_{gi}^{(k)}) d_{p_{gi}}^2 \right)\notag\\
    & \qquad \qquad \qquad + \sum_{(i,j)\in\calL} \frac{1}{2} \bd_{ij}^\top \nabla^2 \calL_{ij}^{(k)} \bd_{ij} \\
    & \text{subject to} \notag \\ 
    & \underline{p}_{gi}^{(k)} \leq d_{p_{gi}} \leq \bar{p}_{gi}^{(k)},\ (g,i)\in\calG, \label{eq:opfqpc1}\\
    & \underline{q}_{gi}^{(k)} \leq d_{q_{gi}} \leq \bar{q}_{gi}^{(k)},\ (g,i)\in\calG, \label{eq:opfqpc2}\\
    & \underline{v}_i^{(k)} \leq d_{w_i} \leq \bar{v}_i^{(k)},\ i\in\calB, \\
    & \underline{\theta}_i^{(k)}\leq d_{\theta_i}\leq \overline{\theta}_i^{(k)},\ i\in\calB,\\
    & \sum_{(g,i)\in\calG} d_{p_{gi}} - \sum_{j\in\calB_i} d_{p_{ij}} - g_i^{sh} d_{w_i} = p_i^{d,(k)},\ i\in\calB, \\
    & \sum_{(g,i)\in\calG} d_{q_{gi}} - \sum_{j\in\calB_i} d_{q_{ij}} + b_i^{sh} d_{w_i} = q_i^{d,(k)},\ i\in\calB, \label{eq:opfqpc6}\\
    % & d_{p_{ij}} - (g_{ii} d_{w_i} + g_{ij} d_{w_{ij}^R} + b_{ij} d_{w_{ij}^I}) = 
    % - p_{ij}^{(k)} + g_{ii} w_i^{(k)} + g_{ij} w_{ij}^{R,(k)} + b_{ij} w_{ij}^{I,(k)} && (i,j) \in \calL \\
    % & d_{q_{ij}} - (-b_{ii} d_{w_i} - b_{ij} d_{w_{ij}^R} + g_{ij} d_{w_{ij}^I} ) = - q_{ij}^{(k)} - b_{ii} w_i^{(k)} - b_{ij} w_{ij}^{R,(k)} + g_{ij} w_{ij}^{I, (k)} && (i,j) \in \calL \\
    % & d_{p_{ji}} - (g_{jj} d_{w_j} + g_{ji} d_{w_{ij}^R} - b_{ji} d_{w_{ij}^I} ) = - p_{ji}^{(k)} + g_{jj} w_j^{(k)} + g_{ji} w_{ij}^{R,(k)} - b_{ji} w_{ij}^{I,(k)} && (i,j) \in \calL \\
    % & d_{q_{ji}} - (-b_{jj} d_{w_j} - b_{ji} d_{w_{ij}^R} - g_{ji} d_{w_{ij}^I} ) = - q_{ji}^{(k)} - b_{jj} w_j^{(k)} - b_{ji} w_{ij}^{R,(k)} - g_{ji} w_{ij}^{I,(k)} && (i,j) \in \calL \label{eq:opfqpc10}\\
    & 2 w_{ij}^{R,(k)} d_{w_{ij}^R} + 2 w_{ij}^{I,(k)} d_{w_{ij}^I} - w_j^{(k)} d_{w_i} - w_i^{(k)} d_{w_j}\notag \\
    &= w_i^{(k)} w_j^{(k)} - (w_{ij}^{R,(k)})^2 - (w_{ij}^{I,(k)})^2 ,\ (i,j)\in\calL, \label{eq:opfqpc11}\\
    & \sin_{ij}^{(k)} d_{w_{ij}^R} - \cos_{ij}^{(k)} d_{w_{ij}^I} \notag\\
    &\quad + (w_{ij}^{R,(k)} \cos_{ij}^{(k)} + w_{ij}^{I,(k)} \sin_{ij}^{(k)}) (d_{\theta_i} - d_{\theta_j}) \notag \\
    &\qquad = -w_{ij}^{R,(k)} \sin_{ij}^{(k)} + w_{ij}^{I,(k)} \cos_{ij}^{(k)}, (i,j)\in\calL, \label{eq:opfqpc12}\\
    & 2 p_{ij}^{(k)} d_{p_{ij}} + 2 q_{ij}^{(k)} d_{q_{ij}} \notag\\
    &\qquad \leq \bar{s}_{ij}^2 - (p_{ij}^{(k)})^2 - (q_{ij}^{(k)})^2 ,\ (i,j)\in \calL, \label{eq:opfqpc13}\\
    & 2 p_{ji}^{(k)} d_{p_{ji}} + 2 q_{ji}^{(k)} d_{q_{ji}} \notag\\
    &\qquad \leq \bar{s}_{ji}^2 - (p_{ji}^{(k)})^2 - (q_{ji}^{(k)})^2 ,\ (i,j)\in\calL, \label{eq:opfqpc14}\\
    & \|\left(d_{p_{gi}}, d_{q_{gi}}, d_{w_i}, d_{\theta_i},d_{w_{ij}^R},d_{w_{ij}^I}\right)\|_\infty \leq \Delta^{(k)},
\end{align}
\end{subequations}
where we denote the vector of independent variables by 
$${\bd}_{ij}:=( d_{w_{i}}, d_{w_{j}}, d_{\theta_{i}}, d_{\theta_{j}})$$
and the variables $(d_{p_{ij}}, d_{q_{ij}}, d_{p_{ji}},d_{q_{ji}})$ corresponding to the power flow variables in \eqref{eq: acopf_pf} are defined as follows.
\begin{subequations}\label{eq: acopf_qppf}
\begin{align}
    & d_{p_{ij}} := g_{ii}(w_i^{(k)} + d_{w_i}) + g_{ij} (w_{ij}^{R,(k)} + d_{w_{ij}^R}) \notag\\
    &\qquad\qquad + b_{ij} (w_{ij}^{I,(k)} + d_{w_{ij}^I}) - p_{ij}^{(k)},\ (i,j) \in \calL,\label{eq: dpij}\\
    & d_{q_{ij}} :=  -b_{ii}(w_i^{(k)} + d_{w_i}) - b_{ij}(w_{ij}^{R,(k)} + d_{w_{ij}^R})\notag\\
    &\qquad\qquad + g_{ij} (w_{ij}^{I,(k)} + d_{w_{ij}^I}) - q_{ij}^{(k)},\ (i,j) \in \calL ,\label{eq: dqij}\\
    & d_{p_{ji}} := g_{jj} (w_j^{(k)} + d_{w_j}) + g_{ji} (w_{ij}^{R,(k)} + d_{w_{ij}^R}) - b_{ji} (w_{ij}^{I,(k)}\notag\\
    &\qquad\qquad + d_{w_{ij}^I}) - p_{ji}^{(k)},\ (i,j) \in \calL, \label{eq: dpji}\\
    & d_{q_{ji}} := -b_{jj} (w_j^{(k)} + d_{w_j})- b_{ji} (w_{ij}^{R,(k)} + d_{w_{ij}^R}) \notag\\
    &\qquad\qquad - g_{ji} (w_{ij}^{I,(k)} + d_{w_{ij}^I}) - q_{ji}^{(k)},\ (i,j) \in \calL \label{eq: dqji}
\end{align}
\end{subequations}
All the other supporting coefficients are defined as follows:
% \begin{subequations}\label{eq: constants}
\begin{align*}
    &\underline{p}_{gi}^{(k)} := \underline{p}_{gi} - p_{gi}^{(k)}, \quad
    \bar{p}_{gi}^{(k)} := \bar{p}_{gi} - p_{gi}^{(k)}, \\
    &\underline{v}_i^{(k)} := \underline{v}_i^2 - w_i^{(k)}, \quad
    \bar{v}_i^{(k)} := \bar{v}_i^2 - w_i^{(k)}, \\
    &\underline{\theta}_i^{(k)} := -2\pi-\theta_i^{(k)},\quad
    \overline{\theta}_i^{(k)} := 2\pi-\theta_i^{(k)},\\
    &p_i^{d,(k)} := p_i^d - \sum_{(g,i)\in\calG} p_{gi}^{(k)} + \sum_{j\in\calB_i} p_{ij}^{(k)} + g_i^{sh} w_i^{(k)}, \\
    &q_i^{d,(k)} := q_i^d - \sum_{(g,i)\in\calG} q_{gi}^{(k)} + \sum_{j\in\calB_i} q_{ij}^{(k)} - b_i^{sh} w_i^{(k)}.
\end{align*}
% \end{subequations}
In practise, $\nabla^2 f_{gi}(p_{gi}^{(k)})>0$ but $\nabla^2\calL_{ij}^{(k)}$ may not be positive semi-definite, hence the QP subproblem \eqref{eq: QPsub} can be nonconvex.
Constraints \eqref{eq:opfqpc1}--\eqref{eq:opfqpc14} and \eqref{eq: acopf_qppf} result from the first-order Taylor approximation of constraints \eqref{eq:opfc1}--\eqref{eq:opfc14} and \eqref{eq: acopf_pf}, respectively.

\subsection{Improved QP Formulation}\label{eq: improv}
We observe that variables $d_{w_{ij}^R}$ and $d_{w_{ij}^I}$ in the QP subproblem \eqref{eq: QPsub} are dependent on the variable vector $\bd$ by the equality constraints \eqref{eq:opfqpc11}--\eqref{eq:opfqpc12}. With that, the trust-region constraints on these dependent variables can lead to QP infeasibility unnecessarily. To avoid the issue, we remove the dependent variables $d_{w_{ij}^R}$ and $d_{w_{ij}^I}$ by substituting by linear equations with respect to ${\bd}_{ij}$, which can be written as
\begin{align}\label{eq: elim}
    d_{w_{ij}^R} = {a_{ij}^R}^\top {\bd}_{ij} + b_{ij}^R\quad \text{and}\quad d_{w_{ij}^I} = {a_{ij}^I}^\top {\bd}_{ij} + b_{ij}^I,
\end{align}
with auxiliary parameters $a_{ij}^R, a_{ij}^I, b_{ij}^R, b_{ij}^I$.
As a result, the trust-region constraints on the removed variables can be eliminated.
Here, we omit the detailed solution representation for simplicity. In Section~\ref{sec: line_kernel}, we will show the detailed decomposition formulation using representation \eqref{eq: elim}. 

% we solve equations \eqref{eq:opfqpc11} and \eqref{eq:opfqpc12} and then represent $d_{w_{ij}^R}$ and $d_{w_{ij}^I}$ as functions of the other variables. This representation enables us to eliminate the trust regions on $d_{w_{ij}^R}$ and $d_{w_{ij}^I}$ because they are guaranteed to be bounded given the trust-region constraints on the other variables. Defining ${\bf \hat{d}}_{ij}=( d_{w_{i}}, d_{w_{j}}, d_{\theta_{i}}, d_{\theta_{j}})$, we denote the solution representation of $d_{w_{ij}^R}$ and $d_{w_{ij}^I}$ as
% \begin{align}\label{eq: elim}
%     d_{w_{ij}^R} = {a_{ij}^R}^\top {\bf \hat{d}}_{ij} + b_{ij}^R\quad \text{and}\quad d_{w_{ij}^I} = {a_{ij}^I}^\top {\bf \hat{d}}_{ij} + b_{ij}^I,
% \end{align}
% with auxiliary parameters $a_{ij}^R, a_{ij}^I, b_{ij}^R, b_{ij}^I$. 

\section{Component-Based Decomposition of the QP Subproblem}\label{sec: ADMM}
We present the application of the ADMM algorithm for solving the QP subproblem based on the decomposition into independent smaller subproblems for each generator, transmission line, and bus. This idea follows the existing ADMM works on component-based decomposition of ACOPF \cite{mhanna2019Adaptive} and its GPU implementation \cite{kim2022accelerated}. The difference is that, instead of decomposing the original ACOPF problem, we decompose the QP subproblem, whose objective function involves a large number of Hessian terms that relate to all the decision variables. In addition, we introduce a different decomposition rule for the QP subproblem. 

\subsection{Coupling Constraints}
% {\color{blue} explain ADMM algorithm and decomposition rule we use for QP}

We follow the idea of component-based decomposition \cite{mhanna2019Adaptive} to decompose the QP subproblem \eqref{eq: QPsub} into independent small problems for each generator, bus, and transmission line. To this end, we introduce the following auxiliary variables and the coupling constraints to the QP subproblem,
\begin{subequations}\label{eq: cpl}
    \begin{align}
        & d_{p_{gi}} = d_{p_{gi(i)}}, \quad
        d_{q_{gi}} = d_{q_{gi(i)}}, \quad (g,i)\in\calG ,\\
        & d_{p_{ij}} = d_{p_{ij(i)}}, \quad
        d_{q_{ij}} = d_{q_{ij(i)}}, \quad (i,j) \in \calL, \\
        & d_{p_{ji}}  = d_{p_{ji(j)}}, \quad
        d_{q_{ji}} = d_{q_{ji(j)}}, \quad (i,j) \in \calL ,\\
        & \bd_{ij} = \bar{\bd}_{ij}\ (i,j) \in \calL,
        % & d_{w_{i(ij)}} = d_{w_i},\ d_{w_{j(ji)}} = d_{w_j},\  (i,j)\in\calL, \\
        % & d_{\theta_{i(ij)}} = d_{\theta_i},\  d_{\theta_{j(ji)}} = d_{\theta_j},\ (i,j)\in\calL. 
    \end{align}
\end{subequations}
where $\bar{\bd}_{ij} := ( d_{w_{i(ij)}}, d_{\theta_{i(ij)}}, d_{w_{j(ji)}}, d_{\theta_{j(ji)}})$.
This allows us to take the augmented Lagrangian of the QP subproblem with respect to the coupling constraints \eqref{eq: cpl}. 
% with the corresponding Lagrangian multiplier $\lambda$ and penalty parameter $\rho>0$. 

\subsection{ADMM Subproblems}
% {\color{blue} how ADMM works and communicate with SQP}

% We split the constraints of \eqref{eq: QPsub} into different classes of the ADMM subproblem. 
% The detailed generator, transmission line, and bus subproblem formulation and their solution approaches are discussed in the following subsections. 
With the auxiliary variables in \eqref{eq: cpl}, the QP subproblem is decomposed into the generator subproblem with variables $(d_{p_{gi}}, d_{q_{gi}})$ for each generator $(g,i)\in\calG$, the line subproblem with variables $(d_{p_{ij}}, d_{q_{ij}}, d_{p_{ji}}, d_{q_{ji}}, \bar{\bd}_{ij})$ for each transmission line $(i,j)\in\calL$, and the bus subproblem with all the other variables in \eqref{eq: cpl} for each bus $i\in\calB$.

\subsubsection{Generator Subproblem}
% {\color{blue} show generator kernel and closed form solution}

For each generator $(g,i)\in\calG$, the generator subproblem is given as the following convex QP with box constraints:
\begin{subequations}\label{eq: gen_kernel}
\begin{align}
    \min_{d_{p_{gi}}, d_{q_{gi}}} \
    & \nabla f_{gi}(p_{gi}^{(k)})d_{p_{gi}} + \frac{1}{2}\nabla^2 f_{gi}(p_{gi}^{(k)}) d_{p_{gi}}^2 + \lambda_{p_{gi}}^{(k)} d_{p_{gi}}\notag \\
    & \ + \lambda_{q_{gi}}^{(k)} d_{q_{gi}} +  \frac{\rho_{p_{gi}}}{2}\|d_{p_{gi}} - d^{(k)}_{p_{gi(i)}}\|^2 \notag\\
    &\ +\frac{\rho_{q_{gi}}}{2}\|d_{q_{gi}} - d^{(k)}_{q_{gi(i)}}\|^2 \\
    \text{s.t.} \quad
    & l_{p_{gi}} \leq d_{p_{gi}} \leq u_{p_{gi}},\ l_{q_{gi}} \leq d_{q_{gi}} \leq u_{q_{gi}},
\end{align}
\end{subequations}
where 
\begin{align*}
    & l_{p_{gi}} := \max\{-\Delta^{(k)}, \underline{p}_{gi}^{(k)}\},
    && u_{p_{gi}} := \min\{\Delta^{(k)}, \overline{p}_{gi}^{(k)}\}, \\
    & l_{q_{gi}} := \max\{-\Delta^{(k)}, \underline{q}_{gi}^{(k)}\},
    && u_{q_{gi}} := \min\{\Delta^{(k)}, \overline{q}_{gi}^{(k)}\}.
\end{align*} 
For the sake of brevity and with a slight abuse of notation, we present the problem \eqref{eq: gen_kernel} in the following generic formulation:
\begin{align}\label{eq: gen_kernel_simp}
    \min_{l\leq x\leq u} \frac{1}{2}x^\top Qx-c^\top x 
\end{align}
Accordingly, the closed-form solution for \eqref{eq: gen_kernel_simp} is obtained by $x^*=\max(l,\min(u,Q^{-1}c))$, where $\max$ and $\min$ are element-wise operations on vectors. 
% \begin{align}
%     d_{p_{gi}}^* &= \max\left(l_{p_{gi}}, \min\left(u_{p_{gi}},
%                     \frac{-\left(\nabla f_{gi}(p_{gi}^{(k)}) + \lambda^{(k)}_{p_{gi}} - \rho_{p_{gi}} d_{p_{gi}(i)}^{(k)} \right)}
%                     {\nabla^2 f_{gi}(p_{gi}^{(k)}) + \rho_{p_{gi}}}\right)\right) \notag\\
%     d_{q_{gi}}^* &= \max\left(l_{q_{gi}}, \min\left(u_{q_{gi}},
%                     \frac{-\left(\lambda^{(k)}_{q_{gi}} - \rho_{q_{gi}} d_{q_{gi}(i)}^{(k)} \right)}
%                     {\rho_{q_{gi}}}\right)\right) \notag 
% \end{align}

% \begin{align}
%     d_{p_{gi}}^* &= \clamp\Bigg(l_{p_{gi}},u_{p_{gi}},\frac{-\left(\nabla f_{gi}(p_{gi}^{(k)}) + \lambda^{(k)}_{p_{gi}} - \rho_{p_{gi}} d_{p_{gi}(i)}^{(k)} \right)}
%                     {\nabla^2 f_{gi}(p_{gi}^{(k)}) + \rho_{p_{gi}}}\Bigg) \notag\\
%     d_{q_{gi}}^* &= \clamp\Bigg(l_{q_{gi}}, u_{q_{gi}},
%                     \frac{-\left(\lambda^{(k)}_{q_{gi}} - \rho_{q_{gi}} d_{q_{gi}(i)}^{(k)} \right)}
%                     {\rho_{q_{gi}}}\Bigg) \notag 
% \end{align}
% where $\clamp\{l,u,x\} := \max\{l,\min\{u,x\}\}$.
\subsubsection{Line Subproblem}\label{sec: line_kernel}

% {\color{blue} show line kernel and explain GPU implementation with augmented lagrangian and Exatron}

% With the component-based decomposition and elimination step in Section~\ref{eq: improv}, 
For each transmission line $(i,j)\in\calL$, the line subproblem is given as the following nonconvex QP subproblem on four variables $\bar{\bd}_{ij}$ with inequality constraints:
\begin{subequations}\label{eq: line_kernel}
\begin{align}
    &\min_{\bar{\bd}_{ij}} \ 
     \frac{1}{2} \bar{\bd}_{ij}^\top A_{ij}^\top H_{ij} A_{ij} \bar{\bd}_{ij} + b_{ij}^\top H_{ij} A_{ij} \bar{\bd}_{ij} \notag\\
      & +  \left[\lambda^{(k)}_{p_{ij}}( {d}_{p_{ij}} - d^{(k)}_{p_{ij(i)}} )+\lambda^{(k)}_{q_{ij}}( {d}_{q_{ij}} - d^{(k)}_{q_{ij(i)}} )\right]\notag\\
        & +  \left[\frac{\rho_{p_{ij}}}{2}\| {d}_{p_{ij}} - d^{(k)}_{p_{ij(i)}}\|^2+\frac{\rho_{q_{ij}}}{2}\| {d}_{q_{ij}} - d^{(k)}_{q_{ij(i)}}  \|^2\right]\notag\\
        & + \left[\lambda^{(k)}_{w_{i_{(ij)}}}(d_{w_{i(ij)}} - d^{(k)}_{w_i})+\lambda^{(k)}_{\theta_{i_{(ij)}}}(d_{\theta_{i(ij)}} - d^{(k)}_{\theta_i}\right]\notag\\
        & + \left[\frac{\rho_{w_{i_{(ij)}}}}{2}\|d_{w_{i(ij)}} - d^{(k)}_{w_i}\|^2+\frac{\rho_{\theta_{i_{(ij)}}}}{2}\|d_{\theta_{i(ij)}} - d^{(k)}_{\theta_i}\|^2\right]\notag\\
        & + \left[\lambda^{(k)}_{p_{ji}}({d}_{p_{ji}}  - d^{(k)}_{p_{ji(j)}})+\lambda^{(k)}_{q_{ji}}({d}_{q_{ji}} - d^{(k)}_{q_{ji(j)}})\right]\notag\\
        & + \left[\frac{\rho_{p_{ji}}}{2}\|{d}_{p_{ji}}  - d^{(k)}_{p_{ji(j)}}\|^2+\frac{\rho_{q_{ji}}}{2}\|{d}_{q_{ji}} - d^{(k)}_{q_{ji(j)}}\|^2\right]\notag\\
        & + \left[\lambda^{(k)}_{w_{j_{(ji)}}}(d_{w_{j(ji)}} - d^{(k)}_{w_j})+\lambda^{(k)}_{\theta_{j_{(ji)}}}(d_{\theta_{j(ji)}} - d^{(k)}_{\theta_j})\right]\notag\\
        & + \left[\frac{\rho_{w_{j_{(ji)}}}}{2}\|d_{w_{j(ji)}} - d^{(k)}_{w_j}\|^2+\frac{\rho_{\theta_{j_{(ji)}}}}{2}\|d_{\theta_{j(ji)}} - d^{(k)}_{\theta_j}\|^2\right]\\
    &\text{s.t.} \
      \underline{v}_i^{(k)} \leq d_{w_{i(ij)}} \leq \bar{v}_i^{(k)},\quad \underline{v}_j^{(k)} \leq d_{w_{j(ji)}} \leq \bar{v}_j^{(k)},\\
    & \underline{\theta}_i^{(k)}\leq d_{\theta_{i(ij)}}\leq \overline{\theta}_i^{(k)},\quad \underline{\theta}_j^{(k)}\leq d_{\theta_{j(ji)}}\leq \overline{\theta}_j^{(k)},\\
    &2 p_{ij}^{(k)} d_{p_{ij}} + 2 q_{ij}^{(k)} d_{q_{ij}} \leq \bar{s}_{ij}^2 - (p_{ij}^{(k)})^2 - (q_{ij}^{(k)})^2,\label{eq: line_kernel_1}\\
    &2 p_{ji}^{(k)} d_{p_{ji}} + 2 q_{ji}^{(k)} d_{q_{ji}} \leq \bar{s}_{ji}^2 - (p_{ji}^{(k)})^2 - (q_{ji}^{(k)})^2,\label{eq: line_kernel_2}\\
     &  \|d_{w_{i(ij)}}\| \leq \Delta^{(k)}, \quad \|d_{\theta_{i(ij)}}\| \leq \Delta^{(k)},\\
    & \|d_{w_{j(ji)}}\| \leq \Delta^{(k)}, \quad \|d_{\theta_{j(ji)}}\| \leq \Delta^{(k)},
\end{align}
\end{subequations}
where $H_{ij}$ represents the corresponding submatrix of the Hessian $\nabla^2 \calL_{ij}^{(k)}$ for line $(i,j)$, and the supporting matrix and vector are constructed as
\begin{align*}
    A_{ij} = \begin{bmatrix}
        \mathbb{I}_2 & & \mathbf{0}_2 \\
        \text{---} & a^R_{ij} & \text{---}\\
         \text{---} & a^I_{ij} & \text{---}\\
        \mathbf{0}_2 & & \mathbb{I}_2 \\
    \end{bmatrix}\in\mathbb{R}^{6\times 4},\quad\quad b_{ij}=\begin{bmatrix}
        \mathbf{0}_2\\b_{ij}^R\\b_{ij}^I\\ \mathbf{0}_2
    \end{bmatrix}\in\mathbb{R}^6,
\end{align*}
where $a_{ij}^R, a_{ij}^I, b_{ij}^R, b_{ij}^I$ are auxiliary parameters from Section~\ref{eq: improv}.
Moreover, by the definitions of \eqref{eq: acopf_qppf} and \eqref{eq: elim}, we substitute $(d_{p_{ij}},d_{q_{ij}},d_{p_{ji}},d_{q_{ji}})$ in the line subproblem \eqref{eq: line_kernel} by
% \begin{subequations}
\begin{align*}
    {d}_{p_{ij}} &= g_{ii}(w_i^{(k)} + d_{w_{i(ij)}}) + g_{ij}(w_{ij}^{R,(k)} + {a_{ij}^R}^\top {\bf \bar{d}}_{ij} + b_{ij}^R)\notag\\
    & \quad + b_{ij} (w_{ij}^{I,(k)} +{a_{ij}^I}^\top {\bf \bar{d}}_{ij} + b_{ij}^I) - p_{ij}^{(k)},  \\
    {d}_{q_{ij}} &= -b_{ii} (w_i^{(k)} + d_{w_{i(ij)}}) - b_{ij} (w_{ij}^{R,(k)} + {a_{ij}^R}^\top {\bf \bar{d}}_{ij} + b_{ij}^R) \notag\\
    & \quad + g_{ij} (w_{ij}^{I,(k)} +{a_{ij}^I}^\top {\bf \bar{d}}_{ij} + b_{ij}^I) - q_{ij}^{(k)},  \\
    {d}_{p_{ji}} &= g_{jj} (w_j^{(k)} + d_{w_{j(ji)}}) + g_{ji} (w_{ij}^{R,(k)} + {a_{ij}^R}^\top {\bf \bar{d}}_{ij} + b_{ij}^R) \notag\\
    & \quad - b_{ji} (w_{ij}^{I,(k)} +{a_{ij}^I}^\top {\bf \bar{d}}_{ij} + b_{ij}^I)  - p_{ji}^{(k)}, \\
    {d}_{q_{ji}} &= -b_{jj} (w_j^{(k)} + d_{w_{j(ji)}}) - b_{ji} (w_{ij}^{R,(k)} + {a_{ij}^R}^\top {\bf \bar{d}}_{ij} + b_{ij}^R) \notag\\
    & \quad - g_{ji} (w_{ij}^{I,(k)} +{a_{ij}^I}^\top {\bf \bar{d}}_{ij} + b_{ij}^I) - q_{ji}^{(k)}. 
\end{align*}
% \end{subequations}
% Note that the line subproblem has only four variables $\bar{\bd}_{ij} = ( d_{w_{i(ij)}}, d_{\theta_{i(ij)}}, d_{w_{j(ji)}}, d_{\theta_{j(ji)}})$ for each $(i,j)\in\calL$.

Unfortunately, no closed-form solution is available for the line subproblem. To solve \eqref{eq: line_kernel}, we use the augmented Lagrangian method (ALM) \cite{ALM1, ALM2, ALM3}. As a result, in each iteration, we solve a box-constrained NLP problem by using a trust-region Newton method (TRON) \cite{lin1999newton}.
TRON has been implemented in Julia for solving a large number of small box-constrained nonlinear programs on GPUs \cite{kim2021Leveraging}. 
To demonstrate the reformulation process, we consider the following generic formulation of \eqref{eq: line_kernel}:
\begin{subequations}\label{eq: line_kernel_simp}
\begin{align}
        \min_{l\leq x \leq u} \quad & \frac{1}{2}x^\top Qx-c^\top x\\
        \text{s.t.} \quad & h(x)\leq 0, \label{eq: line_kernel_simp_2}
\end{align}
\end{subequations}
where \eqref{eq: line_kernel_simp_2} represents the line limit constraints \eqref{eq: line_kernel_1} and \eqref{eq: line_kernel_2}. Introducing auxiliary variable $s$, we rewrite \eqref{eq: line_kernel_simp} as
\begin{subequations}\label{eq: line_kernel_simp1}
\begin{align}
        \min_{l\leq x \leq u, s\geq 0} \quad & \frac{1}{2}x^\top Qx-c^\top x\\
        \text{s.t.} \quad & h(x)+s = 0 \label{eq: line_kernel_simp1_2}.
\end{align}
\end{subequations}
To use ALM, we penalize the equality constraints \eqref{eq: line_kernel_simp1_2} to the objective and obtain 
\begin{align}\label{eq: line_kernel_alm}
        \min_{l\leq x \leq u,
         s\geq 0} \ \frac{1}{2}x^\top Qx-c^\top x -u^\top h(x) + \frac{1}{2\mu}\|h(x)\|_2^2,    
\end{align}
where $u$ and $\mu$ are auxiliary parameters. Then, the boxed-constrained problem \eqref{eq: line_kernel_alm} is solved by TRON. The details of the ALM including the parameter update and termination can be found in \cite{ALM3}.   

% {\color{blue}explain ExaTron and augmented Lagrangian} 

\subsubsection{Bus Subproblem}
% {\color{blue} show bus kernel and close-form solution}

Similarly, for each bus $i\in\calB$, the bus subproblem is given by the following convex QP with only linear equality constraints:
\begin{subequations}\label{eq: bus_kernel}
    \begin{align}
    \min \
    & \sum_{g\in\calG_i} \left[ - \lambda_{p_{gi}}^{(k)} d_{p_{gi(i)}} + \frac{\rho_{p_{gi}}}{2} \left( d_{p_{gi}}^{(k+1)} - d_{p_{gi(i)}} \right)^2 \right] \notag \\
    & + \sum_{g\in\calG_i} \left[ - \lambda_{q_{gi}}^{(k)} d_{q_{gi(i)}} + \frac{\rho_{q_{gi}}}{2} \left( d_{q_{gi}}^{(k+1)} - d_{q_{gi(i)}}\right)^2 \right] \notag \\
    & + \sum_{j\in\calB_i} \left[ - \lambda_{p_{ij}}^{(k)} d_{p_{ij(i)}} + \frac{\rho_{p_{ij}}}{2} \left( d_{p_{ij}}^{(k+1)} -  d_{p_{ij(i)}} \right)^2 \right] \notag \\
    & + \sum_{j\in\calB_i} \left[ - \lambda_{q_{ij}}^{(k)} d_{q_{ij(i)}} + \frac{\rho_{q_{ij}}}{2} \left( d_{q_{ij}}^{(k+1)} -  d_{q_{ij(i)}} \right)^2 \right] \notag \\
    & + \sum_{j\in\calB_i} \left[ - \lambda_{w_{i(ij)}}^{(k)} d_{w_i} + \frac{\rho_{w_{i(ij)}}}{2} \left(  d_{w_{i(ij)}}^{(k+1)} - d_{w_i}\right)^2 \right] \notag \\
    & + \sum_{j\in\calB_i} \left[ - \lambda_{\theta_{i(ij)}}^{(k)} d_{\theta_i} + \frac{\rho_{\theta_{i(ij)}}}{2} \left(   d_{\theta_{i(ij)}}^{(k+1)} - d_{\theta_i} \right)^2 \right] \\
     \text{s.t.} \quad &  \sum_{(g,i)\in\calG} d_{p_{gi(i)}} - g_i^{sh} d_{w_i} = p_i^{d,(k)} + \sum_{j\in\calB_i} d_{p_{ij(i)}},\\
    & \sum_{(g,i)\in\calG} d_{q_{gi(i)}} + b_i^{sh} d_{w_i} = q_i^{d,(k)} + \sum_{j\in\calB_i} d_{q_{ij(i)}}.
\end{align}
\end{subequations}

The Karush--Kuhn--Tucker (KKT) conditions of \eqref{eq: bus_kernel} are linear equality systems that can be solved directly with closed-form solutions. To that end, we rewrite $\eqref{eq: bus_kernel}$ into the following generic form:
\begin{subequations}\label{eq: bus_kernel_simp}
\begin{align}
       \min_{x} \quad &\frac{1}{2}x^\top Qx-c^\top x \\
       \text{s.t.} \quad & Ax = b. \label{eq: bus_kernel_simp_2}
\end{align}
\end{subequations}

By directly solving the KKT conditions of \eqref{eq: bus_kernel_simp}, we get that the optimal multiplier of \eqref{eq: bus_kernel_simp_2} is $\lambda^* = (AQ^{-1}A^\top)^{-1}(AQ^{-1}c-b)$ and the optimal solution $x^*=Q^{-1}(c-A^\top\lambda^*)$.  

\subsubsection{Multiplier Update}
% {\color{blue} how multiplier update}

To introduce the rules of multiplier update and termination in the ADMM algorithm, instead of using the actual coupling constraints \eqref{eq: cpl}, we use a simpler but similar coupling constraint representation $x=\bar{x}$, where $x$ are the decision variables and $\bar{x}$ is the auxiliary variable, with multiplier $\lambda$ and penalty $\rho$. At iteration $k$, we update the multiplier as
\begin{align}\label{eq: multiplier}
    \lambda^{(k+1)} = \lambda^{(k)} + \rho(x^{(k+1)} -\bar{x}^{(k+1)} ),  
\end{align}
where $x^{(k+1)}$, $\bar{x}^{(k+1)}$ can be obtained by solving the decomposed problems defined under $\lambda^{(k)}$,  $x^{(k)}$, $\bar{x}^{(k)}$. In our case they are the generator problem \eqref{eq: gen_kernel}, the transmission line problem \eqref{eq: line_kernel}, and the bus problem \eqref{eq: bus_kernel}.

% {\color{blue} how ADMM terminates}

% To improve the convergence performance, we use a heuristic called residual balancing \cite{admm_residual_balance}. 
\subsubsection{ADMM Algorithm}
% {\color{blue} how ADMM works and communicate with SQP}

With the results in the sections above,  we give the following sketch of the ADMM algorithm. 

To terminate the ADMM, we simply check that the new primal residual $r^{(k+1)}$ and dual residual $s^{(k+1)}$ satisfy the tolerance \cite{ADMM_book}. In our case, $r^{(k+1)} = x^{(k+1)} - \bar{x}^{(k+1)}$ and $s^{(k+1)} = \rho(\bar{x}^{(k+1)} - \bar{x}^{(k)})$. 
After termination, we return the solution of \eqref{eq: QPsub} to the SQP algorithm as well as the current estimate of the multipliers $\pi^{(k)}$ defined in \eqref{eq: lag}. These multipliers can be computed using the KKT conditions of \eqref{eq: line_kernel}.  

\begin{algorithm}[!ht]
\caption{Sketch of the ADMM Algorithm}
\begin{algorithmic}[1]
\REQUIRE Problem parameters of \eqref{eq: QPsub}
\STATE Initialize $d^{(0)}$, $\lambda^{(0)}$, and $\rho>0$ with $k\gets 0$.
\REPEAT
    \STATE Solve the generator problem \eqref{eq: gen_kernel}
    \STATE Solve transmission line problem \eqref{eq: line_kernel}.
    \STATE Solve the bus problem \eqref{eq: bus_kernel}
    \STATE Update multipliers with \eqref{eq: multiplier}.
    \STATE Set $k \gets k+1$.
\UNTIL{termination criteria satisfied}
\end{algorithmic}
\end{algorithm}

% \section{Implementation of the ADMM Algorithm on GPUs}

\section{Numerical Results}\label{sec: result}
We demonstrate the computational performance of our GPU-based SQP algorithm (i.e., SQP-GPU) on different power network cases. We compare with the SQP algorithms with the Ipopt-based QP solver (i.e., SQP-IPOPT) as well as the one with the CPU-based ADMM QP solver (i.e., SQP-CPU). In terms of test cases, we use 15 IEEE power networks on a wide range of network dimension from 9 buses to 6,468 buses \cite{matpower}.

\subsection{Experiment Settings}

We describe the computational test environment as well as the algorithm settings for the SQP algorithm. 
We performed experiments on a workstation equipped with Nvidia’s Quadro GV100 and Intel Xeon 6140 CPU@2.30GHz.

For small instances  (i.e., network dimension $\leq$ 300 buses), we use the SQP penalty $10^5$ and SQP tolerances $10^{-4}$. The remaining ADMM parameter selections will be shown in the result table. 
For larger instances, the SQP-IPOPT method tends to detect infeasible QP subproblems because of the large number of constraints and hence requires a feasibility restoration step. On the other hand, due to the augmented Lagrangian and decomposition nature in ADMM, SQP-ADMM methods can continue to progress regardless of the QP infeasibility. To achieve the best convergence performance, we needed to relax the SQP termination tolerance to $10^{-3}$ or $5\times 10^{-3}$ mainly due to the inaccurate solution of ADMM QP subproblem.
We used $5\times 10^{-3}$ tolerance whenever SQP was not able to reduce the primal infeasibility less than $10^{-3}$ after a number of iterations. 

\begin{table*}[ht]
\centering
\caption{Computation performance comparison of SQP-IPOPT, SQP-CPU, and SQP-GPU.}
\begin{tabular}{
l|rr|rr|rr|rrr|rrr}
\hline
 & \multicolumn{2}{c|}{SQP Steps} & \multicolumn{2}{c|}{SQP Primal Infeas.} & \multicolumn{2}{c|}{SQP Dual Infeas.} & \multicolumn{3}{c|}{ADMM Parameters} & \multicolumn{3}{c}{Time} \\
Case & Ipopt & ADMM & Ipopt & ADMM & Ipopt & ADMM & $\rho$ & Max Iter & $\epsilon$ & Ipopt & ADMM (CPU) & ADMM (GPU) \\
\hline
9   & 3 & 4 & 3.7e-6 & 4.0e-5 & 6.5e-5 & 1.1e-5 & 1e+3 & 1e+3 & 1e-4 & 0.14 & 1.42 & 0.97 \\
30  & 5 & 6 & 2.2e-10 & 1.8e-5 & 2.5e-9 & 7.7e-6 & 2e+4 & 2e+4 & 1e-4 & 0.43 & 29.15 & 11.25 \\
57  & 4 & 3 & 1.8e-8 & 9.4e-5 & 2.3e-9 & 7.2e-5 & 2e+4 & 2e+4 & 1e-4 & 1.51 & 34.03 & 1.85 \\
118 & 3 & 4 & 8.9e-5 & 8.1e-5 & 9.2e-5 & 2.4e-6 & 2e+4 & 1e+3 & 1e-4 & 1.32 & 67.80 & 2.74 \\
300 & 2 & 3 & 1.9e-5 & 7.8e-4 & 7.6e-6 & 1.4e-6 & 2e+4 & 1e+3 & 1e-3 & 2.44 & 111.42 & 2.78 \\
\hline
\end{tabular}\label{tab:small}
\end{table*}

\begin{table*}[ht]
\centering
\caption{Computation performance comparison between IPOPT and SQP-GPU.}
\begin{tabular}{
l|r|rr|rrr|rr}
\hline
 & & \multicolumn{2}{c|}{Infeas.} & \multicolumn{3}{c|}{ADMM Parameters} & \multicolumn{2}{c}{Time} \\
Case & Iters & Primal & Dual & $\rho$ & Max Iter & $\epsilon$ & Ipopt & ADMM \\
\hline
case2383wp  & 2 & 8.7e-4 & 0.0 & 2e+4 & 1e+3 & 1e-3 & 30.70 & 11.79 \\
case2736sp  & 2 & 5.3e-4 & 0.0 & 2e+4 & 1e+3 & 1e-3 & 57.68 & 11.07 \\
case2737sop & 2 & 5.7e-4 & 0.0 & 2e+4 & 1e+3 & 1e-3 & 58.05 & 10.62 \\
case2746wop & 2 & 7.5e-4 & 0.0 & 2e+4 & 1e+3 & 1e-3 & 30.60 & 10.58 \\
case2948rte & 18 & 4.6e-3 & 3.3e-9 & 2e+4 & 1e+3 & 5e-3 & 393.49 & 162.96 \\
case2968rte & 2 & 2.8e-3 & 0.0 & 2e+4 & 1e+3 & 5e-3 & 171.04 & 17.39 \\
case2869pegase & 16 & 3.1e-3 & 7.2e-9 & 2e+4 & 1e+3 & 5e-3 & 61.82 & 125.68 \\
case3012wp  & 3 & 5.6e-4 & 2.0e-7 & 2e+4 & 1e+3 & 1e-3 & 30.78 & 16.02 \\
case3120sp  & 3 & 5.6e-4 & 2.0e-7 & 2e+4 & 1e+3 & 1e-3 & 30.49 & 15.03 \\
case6468rte & 18 & 5.1e-1 & 4.2e-1 & 2e+4 & 1e+3 & 1e-3 & 250.89 & 143.03 \\
\hline
\end{tabular}\label{tab:large}
\end{table*}

% We cannot solve larger instances because SQP generates infeasible QP subproblems and requires a feasibility restoration step.

\subsection{Performance on Small Instances}
% {\color{blue} results table}

The numerical results from five small test cases as well as parameter settings are summarized in Table~\ref{tab:small}. We observe that SQP-IPOPT outperforms the SQP-CPU and SQP-GPU runs due to its high efficiency for solving medium-size QP subproblems. While comparing the SQP-CPU and SQP-GPU, we observe the computation performance advantages resulting from GPU computation. As the system dimension increases, we see that SQP-GPU gains more advantages over SQP-CPU and starts to catch up with SQP-IPOPT.  

% SQP penalty $1e5$; SQP tolerances $1e-4$.

\subsection{Performance on Larger Instances}

For larger instances, the numerical results from ten test cases (up to the system with 6,468 buses) as well as parameter settings are summarized in Table~\ref{tab:large}.  Because of the large network dimension, ADMM algorithm can take more iterations to converge due to the large number of decomposed optimization problems for each component. To address this issue, we relax the termination tolerance and reduce the iteration limit in ADMM to force quick termination, which concurrently drives more iteration steps in SQP. Although early termination in ADMM results in inexact QP solutions and requires more iterations for SQP to converge, we still observe a large performance advantage while comparing with Ipopt\footnote{Here, because SQP-IPOPT will get stuck with QP infeasibility without feasibility restoration, we benchmark the SQP-GPU against a direct solve of ACOPF using Ipopt. Meanwhile, we neglect the results of SQP-CPU as it follows similar patterns as in Table~\ref{tab:small}.}. We see that, in most of cases, SQP-GPU takes much less computational time. Specifically, in case2968rte, we observe an almost 10 times advantages from SQP-GPU over Ipopt. We also observe that the dual infeasibilities are very small even when we are not able to reduce the primal infeasibility below $10^{-3}$.

% {\color{blue}discussion of limitations of restorations and correction to reach high accuracy but preliminary work for the topic}

\subsection{Limitations}

Although we have achieved computational advantages over CPU-based ADMM solvers and Ipopt, there are still limitations to the current SQP-GPU implementation. First, we use only  a basic trust-region SQP structure (i.e., Algorithm~\ref{alg: sqp}) without supporting features such as feasibility restoration and second-order correction \cite{fletcher2002nonlinear}, which may further improve the solution accuracy and computational performance as well as resolving the QP infeasibility. Second, our tuning for ADMM is heuristic and aims to achieve the best performance, which encourages early termination or high tolerance for larger instances. In this sense, we would like to present our work as a preliminary work to this topic of accelerating SQP algorithm using parallel computation.    

% The numerical results from small cases. SQP penalty $1e5$; SQP tolerances $1e-4$.

% For larger instances, we needed to relax the termination tolerance to $10^{-3}$ or $5\times 10^{-3}$ mainly due to the inaccurate solution of ADMM QP subproblem.
% We used $5e-3$ tolerance whenever SQP was not able to reduce the primal infeasibility less than $10^{e-3}$ after a number of iterations. 
% We observe that the dual infeasibilities are very small even when we are not able to reduce the primal infeasibility below $10^{e-3}$.

% We cannot solve larger instances because SQP generates infeasible QP subproblems and requires a feasibility restoration step.

% This additional benefits enables us to avoid the burden of QP infeasibility.

% \footnote{Here, because SQP-IPOPT will get stuck with QP infeasibility, we benchmark against direclty solving ACOPF using Ipopt. Here we also neglect SQP-CPU as it follows similar patterns as in Table~\ref{tab:small}.}. . 

% to $10^{-3}$ or $5\times 10^{-3}$ mainly due to the inaccurate solution of ADMM QP subproblem.
% We used $5e-3$ tolerance whenever SQP was not able to reduce the primal infeasibility less than $10^{-3}$ after a number of iterations. 
% We observe that the dual infeasibilities are very small even when we are not able to reduce the primal infeasibility below $10^{-3}$.

% We cannot solve larger instances because SQP generates infeasible QP subproblems and requires a feasibility restoration step.

% {\color{red}Optional: how tuning affect convergence}

\section{Conclusion}\label{sec: con}
This paper introduces a GPU-accelerated SQP algorithm for solving the ACOPF problem. By combining decomposition and parallel computing techniques, we effectively address the scalability and computational challenges posed by a large quadratic subproblems due to network dimension. Leveraging the ADMM method and the idea of component-based decomposition, we achieve parallel computation of the SQP subproblems on GPU platforms. Through our numerical study, we observed that the GPU-based SQP-ADMM method demonstrates significant reduction in solution time, as compared with the SQP method with the interior-point method QP solver, as well as that with the CPU-based ADMM solver, particularly as the network size increases. As future work, we will focus on techniques to further accelerate the SQP algorithm and improve overall solution accuracy such as addressing the inexact ADMM steps for larger instances. 

% {\color{red}\appendix{optional: Hessian structure}}

% {\color{red}\appendix{optional: Details on improved formulation, e.g., elimination formula}}

% {\color{red}\appendix{optional: close-form solution for bus kernel}}

% The 'triggered' command can be changed if desired:


%\bibliographystyle{IEEEtran}
% \bibitem{Shell}
% M.~Shell, \emph{How to Use the IEEEtran Latex Class}, Latex Archive Contents, \verb+http://www.ieee.org/conferences_events/+ \verb+conferences/publishing/templates.htm+


\bibliographystyle{IEEEtran}
\bibliography{reference.bib}

% \end{thebibliography}

\end{document}